\documentclass{amsart}
\usepackage{graphicx}
\theoremstyle{definition}

\theoremstyle{remark}

\numberwithin{equation}{section}

\begin {document}
\title{ADDITIVE RELATION AND  ALGEBRAIC SYSTEM OF EQUATIONS}
\author{Wu Zi qian}
\address{Xi'an Shiyou University.Address:Xi'an city,China.}
\email{woodschain$@$sohu.com} 
 \pagestyle{myheadings}
\begin{abstract}
\mbox{} Additive relations are defined over additive monoids and additive operation is introduced over these new relations
then we build algebraic system of equations. We can generate profuse equations by additive relations of two variables.
To give an equation with several known parameters is to give an additive relation taking these known parameters as its variables or
value and the solution of the equation is just the reverse of this relation which always exists. We show a core result in this paper
that any additive relation of many variables and their inverse can be expressed in the form of the superposition of additive relations
of one variable in an algebraic system of equations if the system satisfies some conditions. This result means that there is always
a formula solution expressed in the superposition of additive relations of one variable for any equation in this system.
We get algebraic equations if elements of the additive monoid are numbers and get operator equations if they are functions.

\end{abstract}
\maketitle
\section{introduction }

To give explicit solutions are always difficult for not only general algebraic equations but also for general operator
equations. May be one consider that no one try to find formula solution for a general polynomial equation and further
for general algebraic equation after Abel proof that there isn't an algebraic solution for it. But actually the most
outstanding mathematicians devote themselves to this problem in each times.

Camille.Jordan (1838-1922) shows that algebraic equations of any degree can be solved in terms of modular functions in 1870.
Ferdinand.von.Lindemann (1852-1939) expresses the roots of an arbitrary polynomial in terms of theta functions in 1892. In 1895
Emory.McClintock (1840-1916) gives series solutions for all the roots of a polynomial. Robert Hjalmal.Mellin (1854-1933) solves
an arbitrary polynomial equation with Mellin integrals in 1915. In 1925 R.Birkeland shows that the roots of an algebraic equation
can be expressed using hypergeometric functions in several variables. Hiroshi.Umemura expresses the roots of an arbitrary polynomial
through elliptic Siegel functions in1984[3].All of solutions mentioned above are not ones expressed in binary function. David.Hilbert
presumed that there is no solution expressed in binary function for polynomial equations of n when n$\geq$7 and wrote his doubt into
his famous 23 problems as the 13th one[2].

In 1957 V.I.Arnol'd proved that every continuous function can be represented as a superposition of functions of two variables
and refuted Hilbert conjecture[4][5]. Furthermore, A.N.Kolmogorov proved that every continuous function can be represented as a
superposition of continuous functions of one variable and additive operation[1].Thus we can't say we refuse to find formula
solutions for general transcendental equations because they don't exist.

On the hand  modern algebra is running in its own direction but not in the direction of classical mathematics which takes solving equations
including to give formula solutions as one of its main tasks. We can clear this point if we note that so many results gotten by modern
algebra but there are so many algebraic equations with no explicit solution. We will never give explicit solution to most of algebraic
equations  by so few operations which meet axioms of arithmetic and actually  many algebraic equations can't be generated only by them.
But modern algebra limits itself in these operations.

In this paper we limit the problem in finite sets when we construct algebraic equations and operator equations thus the problem becomes
simpler and clearer. We get perfect results then we develop a new algebraic system called  algebraic system of equations by these results.

This is a three grade algebraic system which defines relations on a finite set and defines operations on these relations whereas modern
algebraic structures such as group,ring and field which can be called  two grade algebraic systems are built by set and operations defined
on it. One shall feel rich and colorful of this new algebraic system when he read this paper.

\section{Basic Definitions}

 We shall introduce a new type of relations named  additive relation  and define additive operation on these new relations. We build equations
by these relations and try to solve these equations but we have to face the problem of polykeys in our theory of equation. Function should be the most
suit concept to express polykeys if function can be many-valued but it is defined to be single valued in modern algebra strictly. We don't use function but
use relation to express polykeys then we can avoid conflict with modern algebra. Relations which can be many-valued shall be taken not as operations but
as elements and on these relations we define operations which are single valued thus the new algebraic system will never be inconsistent with modern algebra.

{\bf Definition~2.1} A M+1-ary relation R over non-empty sets $ B_{i}(1\leq i\leq M)$ and $B_{0}$ is a subset of their Cartesian product
written as $R\subset B_{1}\times B_{2}\ldots \times B_{M}\times B_{0}=\{\langle b_{1},b_{2},\ldots,b_{M},b_{0}\rangle|b_{i}\in B_{i}(1\leq i\leq M),b_{0}\in B_{0}\}$.
Specially R is called M+1-ary relation over B if $R\subset \underbrace{B \times B\ldots \times B\times B}_{M+1}=B^{M+1}$ and all M+1-ary relations over B  form the power set of $B^{M+1}$ so we denote them as $P(B^{M+1})$.

{\bf Definition~2.2}  Let B is a finite set and (B,+) is a monoid,R is a M+1-ary relation over B, we shall call it an additive relation of M variables if we take  ith element of its ordered M+1-tuple as its  ith variable and the last element as its value. We denote it by $R^{M}$ and all $R^{M}$ as $P(B^{M+1})$. B is called the basic set of $P(B^{M+1})$ and the number of elements in set B will be taken as the order of $R^{M}$.

The number of all ordered M+1-tuples gotten by B with N elements will be $N^{M+1}$ and the number of $P(B^{M+1})$ will be $2^{(N^{M+1})}$.

For example let $B=\{0,1\}$,then all additive relations of one variables shall be shown below:

$R^{1}_{1}=\emptyset$,$R^{1}_{2}=\{\langle0,0\rangle\}$,$R^{1}_{3}=\{\langle0,1\rangle\}$,$R^{1}_{4}=\{\langle1,0\rangle\}$,$R^{1}_{5}=\{\langle1,1\rangle\}$,
$R^{1}_{6}=\{\langle0,$ $0\rangle,\langle0,1\rangle\}$,$R^{1}_{7}$ $=\{\langle0,0\rangle,\langle1,0\rangle\}$,$R^{1}_{8}=\{\langle0,0\rangle,\langle1,1\rangle\}$,
$R^{1}_{9}=\{\langle0,1\rangle,\langle1,0\rangle\}$, $R^{1}_{10}=\{\langle0,1\rangle,\langle1,1\rangle\}$,$R^{1}_{11}=\{\langle1,0\rangle,\langle1,1\rangle\}$,
$R^{1}_{12}=\{\langle0,0\rangle,\langle0,1\rangle\,\langle1,0\rangle\}$,$R^{1}_{13}=\{\langle0,0\rangle,\langle0,1\rangle,$ $\langle1,1\rangle\}$,
$R^{1}_{14}=\{\langle0,0\rangle,\langle1,0\rangle,\langle1,$$1\rangle\}$,$R^{1}_{15}=\{\langle0,1\rangle,\langle1,0\rangle,\langle1,1\rangle\}$, $R^{1}_{16}=\{\langle0,0\rangle,\langle0,1\rangle,$ $\langle1,0\rangle,\langle1,1\rangle\}$.

$\emptyset$ and $\{\langle0,0\rangle,\langle0,1\rangle,\langle1,0\rangle,\langle1,$ $1\rangle\}$ is also called empty relation and universal relation respectively.

{\bf Definition~2.3} R will be single valued in point $(b_{1},b_{2},\ldots,b_{M})$ if $\exists1\langle b_{1},b_{2},$ $\ldots,b_{M},y\rangle\in R,y\in B$.In this case R will be called function too. R will be many-valued in the same point if $\exists1\langle b_{1},b_{2},\ldots,$ $b_{M},y1\rangle\in R,y1\in B \wedge\exists1\langle b_{1},b_{2},\ldots,b_{M},y2\rangle\in R,y2\in B$. R will be undefined in this point if $\exists0\langle b_{1},b_{2},$ $\ldots,b_{M},y\rangle\in R,y\in B$.

For example,$\{\langle0,0\rangle\}$ and $\{\langle0,1\rangle,\langle1,0\rangle,\langle1,1\rangle\}$ are single valued in point `0'.
$\{\langle0,0\rangle,$ $\langle0,1\rangle$ and $\{\langle0,0\rangle,\langle0,1\rangle,\langle1,0\rangle,\langle1, 1\rangle\}$are many-valued in point `0'.
$\{\langle1,0\rangle\}$ and $\{\langle1,0\rangle,\langle1,1\rangle\}$are undefined in point `0'. $\emptyset$ is undefined in each point.

It's so wonderful to describe the single solution and polykeys and no solution for equations by the single valued and many-valued and undefined of additive relations respectively.

N-th-order additive relations of one variable can be expressed by a table $2\times N$ and elements in the first line are variable  and  ones in the the second line are value.  Many-valued numbers will be partitioned by the symbol `*' and N means undefined.For example $R^{1}_{6}=\{\langle0,0\rangle,\langle0,1\rangle\}$ will be denoted as below. N-th-order additive relations of two variables can be expressed by table(N+1)x(N+1) and we give an example of 3-order one as below too:

\qquad\qquad\qquad\qquad\begin{tabular}{|c|c|}
  \hline
       0&1 \\
  \hline
      0*1&N\\
  \hline
\end{tabular}
{}
\begin{tabular}{|c|c|c|c|}
  \hline
       &0& 1& 2 \\
  \hline
       0&0&1&2\\
  \hline
       1&0*1&0*2&1*2\\
  \hline
       2&$\phi$&0*1*2&$\phi$\\
  \hline
\end{tabular}

For convenience we remaind only the second line and denote $R^{1}_{6}=\{\langle0,0\rangle,\langle0,1\rangle\}$ as (0*1,N).

We use the  expression form for function to denote additive relation of M variables:

\qquad\qquad \qquad\qquad \qquad\qquad$b_{0}=R^{M}(b_{1},b_{2},\ldots , b_{M})$

But to additive relations of two variables which are used so frequently, we sturdy use the form of binary operation $b_{1}R^{2}b_{2}$ to denote it.
It's too clear and we are so accustomed to this form.

{\bf Definition~2.4} Inverse of an additive relation is defined as $R^{-1}=\{\langle x,y\rangle|$ $\langle y,x\rangle\in R\}$.

We denote inverse of R by T(R). Actually inverse means the transformation to elements of additive relation  so we extend it to all transformations of
(1,2,$\ldots$,M,0) and denote transformations as T with superscript $(i_{1},i_{2},\ldots,i_{M},i_{0})$.

$T^{i_{1},i_{2},\ldots,i_{M},i_{0}}(R)=\{\langle b_{i_{1}},b_{i_{2}},\ldots,b_{i_{M}},b_{i_{0}}\rangle | \langle b_{1},b_{2},\ldots ,b_{i},\ldots, b_{M},b_{0} \rangle \in R\}$

There are 6 transformations for additive relations of two variables and  'M!' ones for additive relations of M variables.

Transformations of R can be expressed by the form of function:

\qquad\qquad\qquad\qquad$b_{i_{0}}=T^{i_{1},i_{2},\ldots,i_{M},i_{0}}(R)(b_{i_{1}},b_{i_{2}},\ldots,b_{i_{M}})$

{\bf Definition~2.5}  Composition of additive relations

Let $R^{M}$ is an additive relation of M variables and $\beta$ one of one variable, here we define the ith parameter composition of $R^{M}$
and $\beta$.

For $1\leq i\leq M$ㄩ

$R^{M}\times_{i}\beta=\{\langle b_{1},b_{2},\ldots,b_{i-1},b_{i},b_{i+1},\ldots, b_{M},b_{0}\rangle|y(\langle b_{1},b_{2},$ $\ldots,b_{i-1},y,b_{i+1}, b_{M},$ $b_{0}\rangle\in R^{M}\langle b_{i},y \rangle\in \beta)\}(1\leq i\leq M)$﹝

$R^{M}\times_{i}\beta$ is written in the form of functions:

$[R^{M}\times_{i}\beta](b_{1},b_{2},\ldots,b_{i-1},b_{i},b_{i+1},\ldots, b_{M1})=R^{M}[b_{1},b_{2},\ldots,b_{i-1},\beta( b_{i}),b_{i+1},$ $\ldots, b_{M}]$

For i=0ㄩ

$\beta(b_{0})=R^{M}(b_{1},b_{2},\ldots,b_{i},\ldots, b_{M})$

$[R^{M}\times_{0} \beta](b_{1},b_{2},\ldots,b_{i},\ldots, b_{M})=b_{0}=\beta^{-1}[R^{M}(b_{1},b_{2},\ldots,b_{i},\ldots, b_{M})]$

If both of $R^{M}$ and $\beta$ is an additive relation of one variable  $\beta_{1}$ and  $\beta_{2}$ respectively then:

$b_{0}=[\beta_{1}\times_{0} \beta_{2}](b)=\beta_{2}^{-1}[\beta_{1}(b)]$

{\bf Definition~2.6} Let both of $R^{M}_{1}$ and $R^{M}_{2}$ is additive relation and the sum $R^{M}_{1}+R^{M}_{2}$ is defined as:

$R=(R^{M}_{1}+R^{M}_{2})=\{\langle b_{1},b_{2},\ldots b_{M},b_{01}+b_{02}\rangle|(b_{1},b_{2},\ldots b_{M},b_{01})\in R^{M}_{1} \wedge (b_{1},b_{2},\ldots$  $ b_{M},b_{02})\in R^{M}_{2} \}$

write it in the form of function:

$R=(R^{M}_{1}+R^{M}_{2})(b_{1},b_{2},\ldots,b_{M})=R^{M}_{1}(b_{1},b_{2},\cdots,b_{M})+R^{M}_{2}(b_{1},b_{2},\cdots , b_{M})$

$b_{01}+b_{02}$ must exist and be unique because (B,+) is a monoid. The sum R includes only ordered M+1-tuple so R will exist and be unique despite it may contains
undefined points or many-valued points.

Note,in $P(B^{M+1})$ there is `o' additive relation being `0' value in its any point and o+R=R $(R\in P(B^{M+1}))$. There is domination additive relation being undefined in its any point  and $\emptyset+R=\emptyset (R\in P(B^{M+1}))$. There is local domination additive relation being undefined in some but not all points of it and the sum will be undefined in these points when it adds an other additive relation.

($P(B^{M+1})$,+) is a monoid because there is no an inverse for some additive relations. We can take $P(B^{M+1})$ as a basic set to define additive relations so we can understand why we use not group but monoid in the definition of additive relations.

The sum of two additive relations will be undefined or single valued or many-valued in a point upon the below rule which can be gotten by definition of addition operation.

1 Sum will be undefined in an point if any of additive relations is undefined in this point.

2 Sum will be  single valued in an point if both of additive relations is single valued in this point.

3 Sum will be many-valued in an point if one additive relation is many-valued and the other is defined in this point.

The below example includes all situations:

\begin{tabular}{|c|c|c|c|}
  \hline
       &0& 1& 2 \\
  \hline
       0&$\phi$&$\phi$&$\phi$\\
  \hline
       1&0&1&2\\
  \hline
       2&0*1&1*2&0*1*2\\
  \hline
\end{tabular}
+
{          }
\begin{tabular}{|c|c|c|c|}
  \hline
       &0& 1& 2 \\
  \hline
       0&$\phi$&1&0*1*2\\
  \hline
       1&$\phi$&0&0*2\\
  \hline
       2&$\phi$&0&1*2\\
  \hline
\end{tabular}
     =
{          }
\begin{tabular}{|c|c|c|c|}
  \hline
       &0& 1& 2 \\
  \hline
       0&$\phi$&$\phi$&$\phi$\\
  \hline
       1&$\phi$&1&1*2\\
  \hline
       2&$\phi$&1*2&0*1*2\\
  \hline
\end{tabular}

{\bf Definition~2.7}  Let $R^{M1}$ is an additive relations of M variables and  $R^{M2}$ is a false additive relations of M2 variables defined by $R^{M1}$.

$R^{M2}=F^{M2}_{i_{1},i_{2},\cdots,i_{M1}}(R^{M1})=\{\langle y,\cdots,y,b_{i_{1}},y,\cdots,y,b_{i_{2}},y ,\cdots,y,b_{_{iM1}},$ $y ,\cdots,y,$ $b_{0}\rangle | \langle b_{1},b_{2},\cdots, b_{M1},b_{0} \rangle \in R^{M1}\wedge y\in B,(b_{i_{j}}=b_{j})\}$

Here B is the basic set. $R^{M2}$ will take jth variable of $R^{M1}$ as its $i_{j}$th variable. The value of $R^{M2}$ will change with only M1 variables and other variables are false ones.

For example,$R^{1}=(0,1*2,\emptyset)$, $F^{2}_{1}(R^{1})$ taking variable of  $R^{1}$ as its first variable and $F^{2}_{2}(R^{1})$ taking it as its second variable as below  respectively:

\qquad\qquad\qquad\begin{tabular}{|c|c|c|c|}
  \hline
       &0& 1& 2 \\
  \hline
       0&0&0&0\\
  \hline
       1&1*2&1*2&1*2\\
  \hline
       2&$\emptyset$&$\emptyset$&$\emptyset$\\
  \hline
\end{tabular}
{          }
\begin{tabular}{|c|c|c|c|}
  \hline
       &0& 1&2  \\
  \hline
       0&0&1*2&$\emptyset$\\
  \hline
       1&0&1*2&$\emptyset$\\
  \hline
       2&0&1*2&$\emptyset$\\
  \hline
\end{tabular}

In the expression $w(x,y)=f(x)+g(y)=[F^{2}_{1}(f)+F^{2}_{2}(g)](x,y)$, we change additive relations of one variable `f' and `g' to false one of two variables and
get `w' by adding them then we express `w' by $[F^{2}_{1}(f)+F^{2}_{2}(g)]$ clearly. Certainly we will never get `w' by `f+g'. This denotation is useful when we
express the explicit solution of an equation.

{\bf Definition~2.8} Let B is a finite set and (B,+) is a monoid. We define additive operation `+' over set of additive relations $P(B^{i+1})(1\leq i\leq M)$ defined
over (B,+).[B,$P(B^{i+1})(1\leq i\leq M)$,+] will be called an algebraic system of equations over B.

\section{Basic theorems}

We take additive relations of two variables as examples to show some law below and they are easy to be extended to  ones of many variables.

Here $R^{2}_{i}$ is an additive relation of two variables and $\beta_{i}$ one of one variable.

{\bf Theorem~3.1~~} Addition satisfies commutative law and associative law.

\qquad\qquad\qquad$R^{2}_{1}+R^{2}_{2}=R_{2}^{2}+R^{2}_{1}$

\qquad\qquad\qquad$(R^{2}_{1}+R^{2}_{2})+R^{2}_{3}=R^{2}_{1}+(R^{2}_{2}+R^{2}_{3})$

{\bf Theorem~3.2~~} Composition satisfies the transposal law.

\qquad\qquad\qquad$R^{2}_{1}\times_{i}\beta_{1}\times_{j}\beta_{2}=R^{2}_{1}\times_{j}\beta_{2}\times_{i}\beta_{1}\qquad(i\neq j)$

{\bf Theorem~3.3~~} Addition and composition satisfy left distributive law.

\qquad\qquad\qquad$(R^{2}_{1}+R^{2}_{2})\times_{i}\beta=R^{2}_{1}\times_{i}\beta+R^{2}_{2}\times_{i}\beta\qquad(i\neq 0)$

{\bf Theorem~3.4~~}Transformation satisfies associative law.

\qquad\qquad\qquad$[(c_{1}c_{2})c_{3}](R^{2})=[c_{1}(c_{2}c_{3})](R^{2})$

Transformations make a group with identity element $T^{1,2,0}$ but not an abelian group because it doesn't satisfy the commutative law.

{\bf Theorem~3.5~~} Transformation $T^{2,1,0}$ and composition is equal a new composition and the same transformation.

\qquad\qquad\qquad$[T^{2,1,0}(R^{2})]\times_{1}\beta=T^{2,1,0}(R^{2}\times_{2}\beta)$

\qquad\qquad\qquad$[T^{2,1,0}(R^{2})]\times_{2}\beta=T^{2,1,0}(R^{2}\times_{1}\beta)$

{\bf Theorem~3.6~~} Addition and transformation $T^{2,1,0}$ satisfy left distributive law.

\qquad\qquad\qquad$T^{2,1,0}(R^{2}_{1}+R^{2}_{2})=T^{2,1,0}(R^{2}_{1})+T^{2,1,0}(R^{2}_{2})$

\section{Solvability of an algebraic system of equations}

If an additive relation of M variables can be written in the expression consisting of additive relations of one variable $f_{i},g_{ij}$ like this:

\begin {equation}\label{eq:eps}
R(x_{1},x_{2},\cdots,x_{M})=\sum _{i=1}^{L}f_{i}\big[g_{i1}(x_{1})+g_{i2}(x_{2}) +\cdots+ g_{in}(x_{M})\big]
\end{equation}

then we say that it can be represented as a superposition of additive relations of one variable or decomposing it to form in additive relations of one variable.
We have our core theorem below,a very very important theorem!

{\bf Theorem~4.1~~}$B=\{0,1,\cdots,N\}(N\geq3)$ㄛthen $R^{i}(B)(2\leq i\leq M)$can be represented as a superposition of additive relations of one variable.

{\bf Definition~4.1} A singular additive relation of M variables is one with only one non-zero point.

$\textbf{Proof}$ step1: It holds for the below singular additive relation of two variables.

\qquad \qquad \qquad \qquad \qquad\begin{tabular}{|c|c|c|c|} \hline
& 0& 1 & 2\\
\hline
0&1& 0& 0  \\
\hline
1&0 &0& 0  \\
\hline
2& 0 & 0 & 0 \\
\hline
\end{tabular}

It can be expressed as:

\qquad\qquad \qquad \qquad$R^{2}_{1}(x_{1},x_{2})=(0,0,1)\Big[(1,0,0)(x_{1})+(1,0,0)(x_{2})\Big]$

Step 2: It holds for a general singular additive relation of two variables.

(1,0,0) in (1,0,0)$(x_{1})$ in the expression of $R^{2}_{1}$ is called row additive relation of one variable. Row including the none-zero point will change if we adjust the location of `1' in (1,0,0)$(x_{1})$. (1,0,0) in (1,0,0)$(x_{2})$ in it is called column additive relation of one variable. Column including the none-zero point will change if we adjust the location of `1'in (1,0,0)$(x_{2})$. Both  row additive relation of one variable or  column additive relation of one variable is called location additive relation of one variable.(0,0,1) is called value additive relation of one variable. Value which may be single-valued or many-valued or no-valued will change if we replace `1' with other numbers. Thus we know that it can be represented as a superposition of additive relations of one variable wherever the none-zero point is.

Step 3: It holds for a general 2-th-order additive relation of two variables.

Because every additive relation of two variables can be transformed to sum of 9 singular additive relation of two variables then we get our conclusion.

Step 4: We can extend the result to general N-th-order additive relations of M variables.

In the expression of $R^{2}_{1}$ if we replace row additive relation of one variable (1,0,0)$(x_{1})$ or column additive relation of one variable (1,0,0)$(x_{2})$ by (1,0,$\cdots$,0)$(x_{1})$ or(1,0,$\cdots$,0)$(x_{2})$  and value function (0,0,1) by (0,0,$\cdots$,0,v) respectively then we can extend this expression to order N additive relations of two variables. There will be more location additive relations of one variable (1,0,0, $\cdots$,0) when we extend the result to general N-th-order additive relations of M variables. Here v in   (0,0,$\cdots$,0,v)can be single valued or many-valued or undefined so theorem~4.1 will be hold in any case.

Step 5: If N$<$M+1, the number of location additive relation of one variable in the expression of singular N-th-order additive relation of two variables will be bigger than M+1. For example a 3-th-order singular additive relation of three variables $R^{3}(x_{1},x_{2},x_{3})$ can be represented as:

$R^{3}(x_{1},x_{2},x_{3})=(0,0,1)\Big\{(0,0,1)\big[(1,0,0)(x_{1})+(1,0,0)(x_{2})\big]+(1,0,0)(x_{3})\Big\}$

Note there is an additional location additive relation of one variable (0,0,1) between symbols `$\Big\{$' and `['.

It's easy to check that 2-th-order singular additive relation of two variables can't be represented as a superposition of  additive relations of one variable.

The location additive relation of one variable $g_{ij}$ will not change with W but the value additive relation of one variable $f_{i}$ is not the same for different W so we express $f_{i}$ as below:

\begin {equation}\label{eq:eps}
f_{i}=V_{i}(W)  \qquad   \qquad  \qquad    \qquad (1\leq i\leq K)
\end{equation}

The decomposing method shown above is called trivial method and the number of terms of expression gotten by it is equal to $N^{M}$. The method will be called non-trivial one if the number of terms of expression gotten by it is less than $N^{M}$.

Let B=$\{0,1,2,\cdots N-1\}$, here N is a odd number, additive relations of two variables defined over it can be decomposed to N+1 terms.

$$\sum_{i=0}^{N} h_{i}[f_{i}(x_{1})+g_{i}(x_{2})]=$$

($a_{11}$,$a_{12}$,$a_{13}$,$\cdots$,$a_{1N-1}$,$a_{1N}$)[(1,0,$\cdots$,0,0,0)(x)+(N-2,$\cdots$,3,2,1,0,0)(y)]

+($a_{21}$,$a_{22}$,$\cdots$,$a_{2N-1}$,0)[(0,0,1,2,3,$\cdots$,N-2)(x)+(0,0,$\cdots$,0,0,0,1)(y)]

+($a_{31}$,$a_{32}$,$\cdots$,$a_{3N-1}$,0)[(N-2,$\cdots$,3,2,1,0,0)(x)+(0,0,$\cdots$,0,0,0,1)(y)]

+($a_{41}$,$a_{42}$,$\cdots$,$a_{4N-1}$,0)[(N-1,$\cdots$,4,3,2,1,0)(x)+(1,0,$\cdots$,0,0,0,0)(y)]

+($a_{51}$,$a_{52}$,$\cdots$,$a_{5N-1}$,0)[(N-1,$\cdots$,4,3,2,1,0)(x)+(0,1,$\cdots$,0,0,0,0)(y)]

$\cdots$

+($a_{N1}$,$a_{N2}$,$\cdots$,$a_{N N-1}$,0)[(N-1,$\cdots$,4,3,2,1,0)(x)+(0,0,$\cdots$1,0,0,0)(y)]

+($a_{N+1 1}$,$a_{N+1 2}$,$\cdots$,$a_{N+1 N-1}$,0)[(N-1,$\cdots$,4,3,2,1,0)(x)+(0,0,$\cdots$,0,1,0,0)(y)] (3)

The fourth term and ones before it are called  original items and others plagiarized items. The correctness of the decomposition is easy to be validate by building and solving a set of equations.  Expressions with only N terms may be gotten. It's need to give the procedure to  decompose any additive relation of many variables.

Let $B=\{0,1,\cdots,N\}(N\geq3)$,$A=R^{1}(B)$, [A,$P(A^{i+1})(1\leq i\leq M)$,+] will be an algebraic system of equations over A and there will be many operator equations.

   An algebraic system of equations is solvable if any equation in it can be given a formula solution. Actually  theorem~4.1 give us a conclusion that the system is
solvable. First we judge a new system if it's solvable when we meet it. We shall downplay it if it is unsolvable because there is few elements in the basic set and
is very poor content needs to be studied. There is plenty to be researched in both algebraic system of equations $[B,P(B^{i+1})(1\leq i\leq M),+]$ and $[A,P(A^{i+1})(1\leq i\leq M),+]$. One of both includes so many algebraic equations and another includes so many operator equations. Theorem~4.1  give us a constructive  method to give a formula solution for any equation in that system actually. But the number of terms is too big so it will be a core task for us to get the shortest expression then we can get the perfect formula solution.

We can proof the solvability of $[A,P(A^{i+1})(1\leq i\leq M),+]$ and study equations in it then actually we break a new path for functional analysis.

\section{Formula solution of the double branches equation}

Let $B=\{0,1,2\}$and equation$(xR^{2}_{1}a)R^{2}_{3}$ $(xR^{2}_{2}b)=c$ is called the double branches equation in which $R^{2}_{j} \in P(B^{3})(1\leq j\leq 3)$.
We solve it as the following procedure:

\textbf{Step 1:} Decompose $R^{2}_{3}$ to:
\qquad \qquad\qquad
$$R^{2}_{3}(u,v)=\sum_{i=1}^{4}
f_{i}[g_{i1}(u)+g_{i2}(v)] $$

\qquad \qquad\qquad $$\sum_{i=1}^{4}
f_{i}\Big[g_{i1}(xR^{2}_{1}a)+g_{i2}(xR^{2}_{2}b)\Big]=c$$

\textbf{Step 2:} By composition change $g_{i1}(xR^{2}_{1}a)$ to $x(R^{2}_{1}\times_{0}g_{i1}^{-1})a$ and change $g_{i2}(xR^{2}_{2}b)$ to $x(R^{2}_{2}\times_{0}g_{i2}^{-1})b$:

\qquad\qquad$$\sum_{i=1}^{4}f_{i}\Big[x(R^{2}_{1}\times_{0}g_{i1}^{-1})a+x(R^{2}_{2}\times_{0}g_{i2}^{-1})b\Big]=c$$

\textbf{Step 3:} Change $R^{2}_{1}\times_{0}g_{i1}^{-1}$ and $R^{2}_{2}\times_{0}g_{i2}^{-1}$ to false ones of three variables and add them:

$R^{3}_{3i}(x,a,b)=[F^{3}_{1,2}(R^{2}_{1}\times_{0}g_{i1}^{-1})+F^{3}_{1,3}(R^{2}_{2}\times_{0}g_{i2}^{-1})](x,a,b)\qquad(1\leq i\leq 4)$

\textbf{Step 4:} By composition of $R^{3}_{3i}$ we get:

 \qquad\qquad \qquad
$R^{3}_{4i}(x,a,b)=R^{3}_{3i}\times_{0}f_{i}^{-1}(x,a,b)\qquad(1\leq i\leq 4)$

\textbf{Step 5:} To sum $R^{3}_{4i}$ we get:

\qquad  \qquad  \qquad \qquad \qquad  \qquad
$$R^{3}_{5}(x,a,b)=\sum_{i=1}^{4} R^{2}_{4i}(x,a,b) $$

The equation will become to:

\qquad  \qquad  \qquad \qquad\qquad  \qquad $R^{3}_{5}(x,a,b)=c$

\textbf{Step 6:} By inverse we get:

\qquad \qquad $x=\Big[T^{2,3,0,1}(R^{3}_{5})\Big](a,b,c)=W(a,b,c)$
\\

\textbf{Step 7:} By decomposing the additive relation of three variables we get:

 $$x=\sum_{k=1}^{27}(V_{k}W)\Bigg\{(g_{k4})\Big[(g_{k1})(a)+(g_{k2})(b)\Big]+(g_{k3})(c)\Bigg\}$$

We replace logogram symbols by complete ones.

$$x=\sum_{k=1}^{27}V_{k}\Big[T^{2,3,0,1}\Big(\sum_{i=1}^{4}\big\{\big[F^{3}_{1,2}(R^{2}_{1}\times_{0}g_{i1}^{-1})+F^{3}_{1,3}(R^{2}_{2}\times_{0}g_{i2}^{-1})]\times_{0}(V_{i}R^{2}_{3})^{-1}\big\}\Big)\Big]$$

$$\Bigg\{g_{k4}\Big[T^{2,3,0,1}\Big(\sum_{i=1}^{4}\big\{\big[F^{3}_{1,2}(R^{2}_{1}\times_{0}g_{i1}^{-1})+F^{3}_{1,3}(R^{2}_{2}\times_{0}g_{i2}^{-1})]\times_{0}(V_{i}R^{2}_{3})^{-1}\big\}\Big)\Big]$$

$$\Bigg[g_{k1}\Big[T^{2,3,0,1}\Big(\sum_{i=1}^{4}\big\{\big[F^{3}_{1,2}(R^{2}_{1}\times_{0}g_{i1}^{-1})+F^{3}_{1,3}(R^{2}_{2}\times_{0}g_{i2}^{-1})]\times_{0}(V_{i}R^{2}_{3})^{-1}\big\}\Big)\Big](a)$$

$$+g_{k2}\Big[T^{2,3,0,1}\Big(\sum_{i=1}^{4}\big\{\big[F^{3}_{1,2}(R^{2}_{1}\times_{0}g_{i1}^{-1})+F^{3}_{1,3}(R^{2}_{2}\times_{0}g_{i2}^{-1})]\times_{0}(V_{i}R^{2}_{3})^{-1}\big\}\Big)\Big](b)\Bigg]$$

$$+g_{k3}\Big[T^{2,3,0,1}\Big(\sum_{i=1}^{4}\big\{\big[F^{3}_{1,2}(R^{2}_{1}\times_{0}g_{i1}^{-1})+F^{3}_{1,3}(R^{2}_{2}\times_{0}g_{i2}^{-1})]\times_{0}(V_{i}R^{2}_{3})^{-1}\big\}\Big)\Big](c)\Bigg\}$$

$(xR^{2}_{1}a)R^{2}_{3}$$(xR^{2}_{2}b)=c$ will be an operator equation if we replace the elements of B=$\{0,1,2\}$ by functions of one variable then B=$\{(0,0,0),(2,0,1),$ $(1,0,2)\}$ and we can also give its explicit solution!

We can see that we give the formula solution to this equation not by axioms of arithmetic but by decomposing additive relations of many variables.

\section{Conclusion and expectation}

   We build the algebraic system of equations which is a new  three grade algebraic system  then we can research general additive relations which include all
operations.We show that there is always formula solution for any equation in this system and developed a constructive method for it. Importance  of the results
shown in this paper is twofold.In the first place,the new system study all relations  which include all operations and very few of them  satisfy axioms of arithmetic
and are studied by algebra. So the new system is a great break to algebra. In the second place, we shall find  formula solutions of equations not by axioms of
arithmetic but by expressing additive relations of many variables in the superposition of ones of one variable.

\qquad\qquad \qquad\qquad \qquad\qquad reference

[1] A.N.Kolmogorov, On the representation of continuous functions of
several variables by superpositions of continuous functions of one
variable and addition, Dokl.Akad.Nauk SSSR 114 (1957), 953-956;
English transl., Amer.Math. Soc.Transl. (2) 28 (1963), 55-59.

[2] D.Hilbert, Mathematical Problems, Bull.Amer.Math.
Soc.8(1902),461-462.

[3]H.Umemura. Solution of algebraic equations in terms of theta
constants. In D.Mumford, Tata.Lectures on Theta II, Progress in
Mathematics 43, Birkh user, Boston, 1984.

[4] V.I.Arnol＊d, On functions of three variables, Dokl.Akad. Nauk
SSSR 114 (1957), 679每681; English transl., Amer. Math. Soc.
Transl.(2) 28 (1963), 51每54.

[5] V.I.Arnol＊d, On the representation of continuous functions of
three variables by superpositions of continuous functions of two
variables, Mat.Sb. 48 (1959), 3每74; English transl., Amer.Math.
Soc.Transl.(2) 28 (1963), 61每147.

\end {document}